\documentclass{article}

\usepackage{graphicx}
\usepackage{amssymb}
\usepackage{amsthm}
\usepackage{amsfonts}
\usepackage{marvosym}

\input xy
\xyoption{all}

\usepackage{times}

\usepackage{amsmath}
\usepackage{hyperref}

\newfont{\msbm}{msbm10 scaled\magstephalf}
\def\LS{\rm LS}

\def\cf{\mathop{\rm cf}}
\newtheorem{theorem}{Theorem}[subsection]

\newtheorem{cor}[theorem]{Corollary}

\newtheorem{lemma}[theorem]{Lemma}
\newtheorem{definition}[theorem]{Definition}
\newtheorem{claim}[theorem]{Claim}

\newtheorem{fact}[theorem]{Fact}
\newtheorem{remark}[theorem]{Remark}
\newtheorem{question}[theorem]{Question}

\newtheorem{notation}[theorem]{Notation}

\newtheorem{prop}[theorem]{Proposition}

\newcommand{\seq}[1]{\langle #1 \rangle}

\def\bar#1{\overline{#1}}
\def\implies{\rightarrow}

\def\b1K{\mbox{\boldmath $K$}_{-1}}
\def\bK{\mbox{\boldmath $K$}}

\def\Pscr{{\cal P }}

\def\Kscr{{\cal K }}
 
  \def\Mod{\mathop{\rm Mod}}

\def\bc{\mbox{\boldmath $c$}}

\def\bN{\mbox{\boldmath $N$}}

\def\rest{\mid}

\newbox\noforkbox \newdimen\forklinewidth
\forklinewidth=0.3pt \setbox0\hbox{$\textstyle\smile$}
\setbox1\hbox to \wd0{\hfil\vrule width \forklinewidth depth-2pt
 height 10pt \hfill
2pt\box0\relax}

\def\nonfor\bK_#1{\unionstic\bK_{\textstyle #1}}
\setbox0\hbox{$\textstyle\smile$} \setbox1\hbox to \wd0{\hfil{\sl
/\/}\hfil} \setbox2\hbox to \wd0{\hfil\vrule height 10pt depth
-2pt width \forklinewidth\hfil} \wd1=0 cm \wd2=0 cm
\newbox\doesforkbox
\setbox\doesforkbox\hbox{\lower 2pt\box1 \lower
2pt\box2\lower2pt\box0\relax}

\def\for\bK_#1{\nunionstic\bK_{\textstyle #1}}

\def\LS{\rm LS}

\def\subm{\prec_{\bK}}
\def\grpf #1 #2{{\rm grp}_{#2}(#1)}
\def\spanf #1 #2{{\rm span}_{#2}(#1)}
\def\fldf #1 #2{{\rm fld}_{#2}(#1)}
\def\dclf #1 #2{{\rm dcl}_{#2}(#1)}
\def\rclf #1 #2{{\rm rcl}_{#2}(#1)}
\def\aclf #1 #2{{\rm acl}_{#2}(#1)}
\def\acff #1 #2{{\rm acf}_{#2}(#1)}
\def\strf #1 #2{{\rm str}_{#2}(#1)}
\def\tclf #1 #2{{\rm acf}_{#2}(#1)}

\newcommand{\bb}{\bold{b}}
\newcommand{\ba}{\bold{a}}

\def\hbar{{\bf h}}

\def\xbar{{\bf x}}

 \def\restrict{|}

\newcommand{\sidebar}[1]{\vskip5pt\noindent
 \hskip.50truein\vrule width2.0pt\hskip.5em
 \vbox{\hsize= 4truein\noindent\footnotesize\relax #1 }\vskip5pt\noindent}

\newcommand{\bm}{\bold{m}}
\newcommand{\bn}{\bold{n}}
\newcommand{\bx}{\bold{x}}
\newcommand{\by}{\bold{y}}
\newcommand{\comment}[1]{}

\newcommand{\mcolon}{\!:\!}

\title{Hanf Numbers and Presentation Theorems in AECs}

\author{John Baldwin\thanks{Research partially supported by Simons
travel grant G5402.} \\
Department of Mathematics, Statistics and  Computer
Science\\
University of Illinois at Chicago\\{}\\
Will Boney\thanks{This material is based
 upon work supported by the National Science
  Foundation under Grant No. DMS-1402191. }\\
  Mathematics Department\\
Harvard University}

\begin{document}
\maketitle

\section{Introduction}

This paper addresses a number of fundamental problems in logic and
the philosophy of mathematics by considering some more technical
problems in model theory and set theory.  The interplay between
syntax and semantics is usually considered the hallmark of model
theory. At first sight, Shelah's notion of abstract elementary class
shatters that icon.  As in the beginnings of the modern theory of
structures (\cite{Corry}) Shelah studies certain classes of models
and relations among them, providing an axiomatization in the Bourbaki
(\cite{Bourbaki}) as opposed to the G{\"o}del or Tarski sense:
mathematical requirements, not sentences in a formal language. This
formalism-free approach (\cite{Kennedyff}) was designed to circumvent
confusion arising from the syntactical schemes of infinitary logic;
if a logic is closed under infinite conjunctions, what is the sense
of studying types? However, Shelah's presentation theorem and more strongly Boney's use \cite{Boneyct} of aec's as theories
of $L_{\kappa,\omega}$ (for $\kappa$ strongly compact) reintroduce
syntactical arguments. The issues addressed in this paper trace to
the failure of infinitary logics to satisfy the {\em upward}
L{\"o}wenheim-Skolem theorem or more specifically the compactness
theorem.  The compactness theorem allows such basic algebraic notions
as amalgamation and joint embedding to be easily encoded in first
order logic.  Thus, all complete first order theories have
amalgamation and joint embedding in all cardinalities. In contrast
these and other familiar concepts from algebra and model theory turn
out to be heavily cardinal-dependent for infinitary logic and
specifically for abstract elementary classes. This is especially
striking as one of the most important contributions of modern model
theory is the freeing of first order model theory from its
entanglement with axiomatic set theory (\cite{Baldwindiv}, chapter 7
of \cite{Baldwinphilbook}).

Two main issues are addressed here. We consider not the interaction
of syntax and semantics  in the usual formal language/structure
dichotomy but methodologically. What are reasons for adopting
syntactic and/or semantic approaches to a particular topic? We
compare methods from the very beginnings of model theory with
semantic methods powered by large cardinal hypotheses. Secondly, what
then are the connections of large cardinal axioms with the cardinal
dependence of algebraic properties in model theory. Here we describe
the opening of the gates for  potentially large interactions between
set theorists (and incidentally graph theorists) and model theorists.
More precisely, can the combinatorial properties of small large
cardinals be coded as structural properties of abstract elementary
classes so as to produce Hanf numbers intermediate in cardinality
between `well below the first inaccessible' and `strongly compact'?

Most theorems in mathematics are either true in a specific small
cardinality (at most the continuum) or in all cardinals. For example
 all, {\em finite} division rings are
commutative, thus all finite Desarguesian planes are Pappian. But
{\em all} Pappian planes are Desarguean and not conversely. Of course
this stricture does not apply to set theory, but the distinctions
arising in set theory are combinatorial. First order model theory, to
some extent, and Abstract Elementary Classes (AEC) are beginning to
provide a deeper exploration of Cantor's paradise: algebraic
properties that are cardinality dependent. In this article, we
explore whether certain key properties (amalgamation, joint
embedding, and their relatives) follow this line.
 These algebraic properties are structural in the sense of \cite{Corrybook}.

 Much of this issue arises from an interesting decision of Shelah.  Generalizing
  Fra{\" \i}ss\'{e} \cite{Fraisse} who considered only finite and
  countable stuctures,
 J{\'o}nsson
 laid the foundations for AEC by his study of universal and homogeneous relation systems
  \cite{Jonssonuniv,Jonssonhomuniv}.  Both of these authors assumed
  the amalgamation property (AP) and the joint embedding property (JEP), which in their context is cardinal
  independent.  Variants such as disjoint or free amalgamation (DAP) are a
  well-studied notion in model theory and universal algebra.  But
  Shelah omitted the requirement of amalgamation in defining AEC.
  Two reasons are evident for this: it is cardinal dependent in this
  context; Shelah's theorem (under weak diamond) that categoricity in
  $\kappa$ and few models in $\kappa^+$ implies amalgamation in
  $\kappa$ suggests that amalgamation might be a dividing line.

Grossberg \cite[Conjecture 9.3]{GrossbergBilgi} first raised the
question of the existence of Hanf numbers for joint embedding and
amalgamation in Abstract Elementary Classes (AEC).
 We define four kinds of amalgamation properties (with various cardinal
parameters)
 in
Subsection~\ref{prelim} and a fifth at the end of Section~\ref{syn1}.
The first three notions are staples of the model theory and universal
algebra since the fifties and treated for first order logic in a
fairly uniform manner by the methods of Abraham Robinson. It is a
rather striking feature of Shelah's presentation theorem that issues
of disjointness require careful study for AEC, while disjoint
amalgamation is trivial for complete first order theories.

 Our main result is the following:

\begin{theorem}\label{mt}  Let $\kappa$ be strongly compact and $\bK$ be an
 AEC with L\"{o}wenheim-Skolem
number
 less than $\kappa$.
If $\bK$ satisfies\footnote{This alphabet soup is decoded in
Definition \ref{mdef}.} AP/JEP/DAP/DJEP/NDJEP for models of size
$[\mu,<\kappa)$, then $\bK$ satisfies
    AP/JEP/DAP/DJEP/NDJEP for all models of size $\geq \mu$.

\end{theorem}

 We conclude with a survey of results
showing the large gap for many properties between the largest
cardinal where an `exotic' structure exists and the smallest where
eventual behavior is determined.  Then we provide specific question
to investigate this distinction.

Our starting place for this investigation was second author's work \cite{Boneyct} that emphasized the role of large
cardinals in the study of AEC.  A key aspect of the definition of AEC is as a mathematical
definition with no formal syntax - class of structures satisfying
certain closure properties.  However, Shelah's Presentation Theorem says that AECs are expressible in infinitary languages, $L_{\kappa,\omega}$, which
allowed a proof via sufficiently complete ultraproducts that, assuming enough strongly compact cardinals, all
AEC's were eventually tame in the sense of \cite{GrVaD}.

Thus we approached the problem of finding a Hanf number for
amalgamation, etc.\! from two directions: using ultraproducts to give
purely semantic arguments and using Shelah's Presentation Theorem to
give purely syntactic arguments.  However, there was a gap: although
syntactic arguments gave characterizations similar to those found in
first order, they required looking at the \emph{disjoint} versions of
properties, while the semantic arguments did not see this difference.

The requirement of disjointness in the syntactic arguments stems from a lack of canonicity in Shelah's Presentation Theorem: a single model has many expansions which means that the transfer of structural properties  between an AEC $\bK$ and it's expansion can break down.  To fix this problem, we developed a new presentation theorem, called the \emph{relational presentation theorem} because the expansion consists of relations rather than the Skolem-like functions from Shelah's Presentation Theorem.
\begin{theorem}[The relational presentation theorem, Theorem \ref{preb}] To each AEC
$\bK$ with $LS(K) = \kappa$ in vocabulary $\tau$,  there is an
expansion of $\tau$ by predicates of arity  $\kappa$ and
a theory $T^*$ in  $\mathbb{L}_{(2^\kappa)^+, \kappa^+}$ such that
$\bK$ is exactly the class of $\tau$ reducts of models of $T^*$.
\end{theorem}

Note that this presentation theorem works in $\mathbb{L}_{(2^\kappa)^+, \kappa^+}$ and has symbols of arity $\kappa$, a far cry from the $\mathbb{L}_{(2^\kappa)^+, \omega}$ and finitary language of Shelah's Presentation Theorem.  The benefit of this is that the expansion is canonical or functorial (see Definition \ref{funcexp}).  This functoriality makes the transfer of properties between $\bK$ and $(\Mod T^*, \subset_{\tau^*})$ trivial (see Proposition \ref{abstracttrans}).  This allows us to formulate natural syntactic conditions for our structural properties.

Comparing the relational presentation theorem to Shelah's, another well-known advantage of Shelah's is that it allows for the computation of Hanf numbers for existence (see Section \ref{biggap-sec}) because these exist in $\mathbb{L}_{\kappa,\omega}$.  However, there is an advantage of the relational presentation theorem
: Shelah's Presentation Theorem works with a sentence in the logic $\mathbb{L}_{(2^{LS(\bK)})^+, \omega}$ and there is little hope of bringing that cardinal down\footnote{Indeed an AEC $\bK$ where the sentence is in a smaller logic would likely have to have satisfy the very strong property that there are $<2^{LS(\bK)}$ many $\tau(\bK)$ structures that {\bf are not} in $\bK$}.  On the other hand, the logic and size of theory in the relational presentation theorem can be brought down by putting structure assumptions on the class $\bK$, primarily on the number of nonisomorphic extensions of size $LS(\bK)$, $|\{ (M, N)/\cong : M \prec_{\bK} N \text{ from }\bK_{LS(\bK)}\}|$.  

  We would like to thank Spencer Unger and Sebastien Vasey for helpful discussions regarding these results.

\subsection{Preliminaries}\label{prelim}

We discuss the relevant background of AECs, especially for the case
of disjoint amalgamation.

\begin{definition}\label{mdef} We consider several variations on the joint
embedding property, written JEP or JEP$[\mu, \kappa)$ .
\begin{enumerate}

	\item Given a class of cardinals $\mathcal{F}$ and an AEC $\bK$, $\bK_{\mathcal{F}}$ denotes the collection of $M \in \bK$ such that $|M| \in \mathcal{F}$.  When $\mathcal{F}$ is a singleton, we write $\bK_{\kappa}$ instead of $\bK_{\{\kappa\}}$.  Similarly, when $\mathcal{F}$ is an interval, we write $<\kappa$ in place of $[\LS(\bK), \kappa)$; $\leq \kappa$ in place of $[\LS(\bK), \kappa]$; $>\kappa$ in place of $\{\lambda \mid \lambda > \kappa\}$; and $\geq \kappa$ in place of $\{\lambda \mid \lambda \geq \kappa\}$.

\item An AEC $(\bK,\subm)$ has the joint embedding property, JEP,
    (on the interval $[\mu, \kappa)$)   if any two models (from $\bK_{[\mu, \kappa)}$) can be
    $\bK$-embedded into a larger model.

\item If the embeddings witnessing the joint embedding property can be chosen to have disjoint ranges, then we call this the {\em disjoint
    embedding property} and write $DJEP$.

\item An AEC $(\bK,\subm)$ has the amalgamation property, AP,
    (on the interval $[\mu, \kappa)$)   if, given any triple of models $M_0 \prec M_1, M_2$ (from $\bK_{[\mu, \kappa)}$), $M_1$ and $M_2$ can be
    $\bK$-embedded into a larger model by embeddings that agree on $M_0$.

\item If the embeddings witnessing the amalgamation property can
    be chosen to have disjoint ranges except for $M_0$, then we
    call this the {\em disjoint amalgamation property} and write
    $DAP$.

\end{enumerate}
\end{definition}

\begin{definition}\label{ECG} \begin{enumerate}
\item A {\em finite diagram} \index{finite diagram}or
 {\em $EC(T,\Gamma)$-class} \index{$EC(T,\Gamma)$-class} is the class
 of
models of a first order theory $T$ which omit all types from a
specified collection $\Gamma$ of complete types in finitely many
variables over the empty set.

\item Let $\Gamma$ be  a collection of first order types in finitely
    many variables over the empty set  for a first order theory $T$ in a vocabulary $\tau_1$. A {\em
    $PC(T,\Gamma, \tau)$ class} is the class
 of reducts to $\tau \subset \tau_1$ of
models of a first order $\tau_1$-theory $T$ which omit all
members of  the specified collection $\Gamma$ of partial types.

\end{enumerate}
\end{definition}

\section{Semantic arguments}\label{semanticarg}

It turns out that the Hanf number computation for  the amalgamation
properties is immediate from Boney's ``{\L}o\'s' Theorem for AECs''
\cite[Theorem 4.3]{Boneyct}. We will sketch the argument for
completeness.  For convenience here, we take the following of the
many equivalent definitions of strongly compact; it is the most
useful for ultraproduct constructions.

\begin{definition}[\cite{Jechst3rd}.20] The cardinal $\kappa$ is strongly compact iff for every $S$ and
every $\kappa$-complete filter on $S$ can be extended to a
$\kappa$-complete ultrafilter.  Equivalently, for every $\lambda \geq
\kappa$, there is a fine\footnote{$U$ is fine iff $G(\alpha) := \{z
\in P_\kappa(\lambda) | \alpha \in  z\}$ is an element of $U$  for
each $\alpha < \lambda$.}, $\kappa$-complete ultrafilter on $P_\kappa
\lambda = \{\sigma \subset \lambda : |\sigma|<\kappa\}$.
\end{definition}

For this paper, ``essentially below $\kappa$'' means ``$LS(K) <
\kappa$.''

\begin{fact}[\L o\'{s}' Theorem for AECs] \label{losaec}
Suppose $K$ is an AEC essentially below $\kappa$ and $U$ is a
$\kappa$-complete ultrafilter on $I$. Then $K$ and the class of
$K$-embeddings are closed under $\kappa$-complete ultraproducts and the
ultrapower embedding is a $K$-embedding.
\end{fact}

The argument for Theorem \ref{mt1} has two main steps. First, use
Shelah's presentation theorem to interpret the AEC into
$L_{\kappa,\omega}$ and then use the fact that $L_{\kappa,\omega}$
classes are closed under ultraproduct by {\em $\kappa$-complete}
ultraproducts.

\begin{theorem}\label{mt1}
Let $\kappa$ be strongly compact and $\bK$ be an AEC with
L\"{o}wenheim-Skolem number
 less than $\kappa$.
 \begin{itemize}
\item If $\bK$ satisfies $AP(<\kappa)$ then $\bK$ satisfies $AP$.
	\item If $\bK$ satisfies $JEP(<\kappa)$ then $\bK$ satisfies $JEP$. 	
	\item If $\bK$ satisfies $DAP(<\kappa)$ then $\bK$ satisfies
$DAP$.

\end{itemize}
\end{theorem}

{\bf Proof:}  We first sketch the proof for the first item, $AP$, and
then note the modifications for the other two.

 Suppose that $\bK$ satisfies $AP(<\kappa)$ and
consider a triple of models $(M,M_1,M_2)$ with $M \subm M_1,M_2$ and
$|M|\leq |M_1| \leq |M_2| =\lambda \geq \kappa$.
Now we will use our strongly compact cardinal.
An {\em approximation} of $(M,M_1,M_2)$ is a triple
$\bN =
(N^{\bN}, N_1^{\bN}, N_2^{\bN}) \in (K_{< \kappa})^3$ such that $N^{\bN} \prec M,
N_\ell^{\bN} \prec M_\ell, N^{\bN} \prec N_\ell^{\bN}$ \text{ for } $\ell
= 1, 2$.
We will take an
ultraproduct indexed by the set $X$ below of approximations to the
triple $(M,M_1,M_2)$.
Set
$$X := \{ \bN  \in (\bK_{< \kappa})^3 :  \bN\text{ is an approximation of }(M, M_1, M_2)\}$$

For each $\bN \in X$, $AP(<\kappa)$ implies there is an amalgam of
this triple.  Fix $f^{\bN}_\ell: N^{\bN}_\ell \to N_*^{\bN}$ to
witness this fact.  For each $(A, B, C) \in [M]^{<\kappa} \times
[M_1]^{< \kappa} \times [M_2]^{<\kappa}$, define
$$G(A, B, C):= \{ {\bN} \in X:A \subset N^{\bN}, B \subset N_1^{\bN}, C \subset N^{\bN}_2\}$$
These sets generate a $\kappa$-complete filter on $X$, so it can be
extended to a $\kappa$-complete ultrafilter $U$ on $X $; note that this ultrafilter will satisfy the appropriate generalization of fineness, namely that $G(A,B,C)$ is always a $U$-large set.

We will now take the ultraproduct of the approximations and their amalgam.  In the end, we will end up with the following commuting diagram, which provides the amalgam of the original triple.

\[
\xymatrix{ & M_1 \ar[rr]^{h_1} &  & \Pi N_1^{\bN}/U \ar[dr]^{\Pi f_1^{\bN}} & \\
M \ar[ur] \ar[dr] \ar[rr]_{h} &  & \Pi N^{\bN}/U \ar[ur] \ar[dr] &  & \Pi N_*^{\bN}/U\\
 & M_2 \ar[rr]_{h_2} &  & \Pi N_2^{\bN}/U \ar[ur]_{\Pi f_2^{\bN}} & }
\]

First, we use \L o\'{s}' Theorem for AECs to get the following maps:
\begin{eqnarray*}
h:M \to \Pi N^{\bN}/U & &\\
h_\ell:M_\ell \to \Pi N_\ell^{\bN}/U & & \text{ for }\ell = 1, 2
\end{eqnarray*}
$h$ is defined by taking $m \in M$ to the equivalence class of constant function $\bN \mapsto x$; this constant function is not always defined, but the fineness-like condition guarantees that it is defined on a $U$-large set (and $h_1, h_2$ are defined similarly).  The uniform definition of these maps imply that $h_1 \rest M = h \rest M = h_2 \rest M$.

Second, we can average the $f^{\bN}_\ell$
maps to get ultraproduct maps
$$\Pi f^{\bN}_\ell: \Pi N_\ell^{\bN} / U \to \Pi N_*^{\bN} / U$$
These maps agree on $\Pi N^{\bN}/U$ since each of the individual functions do.  As each $M_\ell$ embeds in $\Pi N_\ell^{\bN} / U$ the composition of the $f$ and $h$ maps gives the amalgam.

There is no difficulty if one of $M_0$ or $M_1$ has cardinality
$<\kappa$; many of the approximating triples will have the same first
or second coordinates but this causes no harm.  Similary, we get the
JEP transfer if $M_0 = \emptyset$. And we can transfer disjoint
amalgamation since in that case each  $N_1^{\bN} \cap N_2^{\bN} =
N^{\bN}$ and this is preserved by the ultraproduct.
$\dag_{\ref{mt1}}$

\section{Syntactic Approaches}

The two methods discussed in this section both depend on expanding
the models of $\bK$ to models in a larger vocabulary.  We begin with
a concept  introduced in Vasey \cite[Definition 3.1]{vasinfstab}.

\begin{definition}\label{funcexp} A {\em functorial expansion} of an AEC $\bK$ in a vocabulary $\tau$
is an AEC $\hat \bK$ in a vocabulary $\hat \tau$ extending $\tau$
such that \begin{enumerate}\item each $M\in \bK$ has a {\em unique}
expansion to a $\hat M \in \hat \bK$, \item  if  $f: M\cong M'$
then $f:\hat M\cong \hat M'$, and \item if $M$ is  a strong
substructure of $M'$ for $\bK$, then $\hat M$ is strong substructure
of  $\hat M'$ for $\hat \bK$.
\end{enumerate}
\end{definition}

This concept unifies a number of previous expansions: Morley's adding
a predicate for each first order definable set, Chang adding a
predicate for each $L_{\omega_1,\omega}$ definable set, $T^{eq}$,
\cite{CHL} adding predicates $R_n(\xbar,y)$ for closure (in an
ambient geometry) of $\xbar$, and the expansion by naming the orbits
in Fra\^{\i}ss\`{e} model\footnote{This has been done for years but
there is a slight wrinkle in e.g. \cite{BKL} where the orbits are not
first order definable.}.

An important point in both \cite{vasinfstab} and our relational
presentation is that the process does not just reduce the complexity
of already definable sets (as Morley, Chang) but adds new definable
sets.  But the crucial distinction here is that the expansion in
Shelah's presentation theorem is not `functorial' in the sense here: each model has several expansions, rather than a single expansion.  That is
why there is an extended proof for amalgamation transfer in
Section~\ref{syn1}, while the transfer in  Section~\ref{relpres}
follows from the following result which is easily proved by chasing
arrows.

\begin{prop}\label{abstracttrans} Let $\bK$ to $\hat \bK$ be a functorial expansion.
$(K, \prec)$ has $\lambda$-amalgamation [joint embedding, etc.] iff
$\hat \bK$ has $\lambda$-amalgamation [joint embedding,
etc.].
\end{prop}

\subsection{Shelah's Presentation Theorem}\label{syn1}

In this section, we provide syntactic characterizations of the
various amalgamation properties in a finitary language.  Our first approach to these results stemmed from the realization that the amalgamation property has the same
syntactic characterization for $L_{\kappa,\kappa}$ as for first order logic if $\kappa$ is strongly compact, i.e.,
the compactness theorem hold for $L_{\kappa,\kappa}$. Combined with Boney's recognition that one could code each AEC with
L{\"o}wenheim-Skolem number less than $\kappa$ in $L_{\kappa,\kappa}$ this seemed a path to showing amalgamation.  Unfortunately,
this path leads through the trichotomy in Fact~\ref{pres}.
 The results
depend directly (or with minor variations) on Shelah's Presentation
Theorem and illustrate its advantages (finitary language) and
disadvantage (lack of canonicity).

\begin{fact}[Shelah's presentation theorem]\label{pres} If $\bK$ is  an AEC (in a  vocabulary $\tau$ with $|\tau| \leq \LS(\bK)$)
 with L\"owenheim-Skolem number $\LS(\bK)$, there is a vocabulary $\tau_1\supseteq \tau$ with cardinality
$|\LS(\bK)|$, a first order $\tau_1$-theory $T_1$ and a set $\Gamma$
of at most $2^{\LS(\bK)}$ partial types such that
\begin{enumerate}

	\item $\bK =\{M'\restrict \tau \mcolon  M' \models T_1 {\ \rm and \ } M'
{\ \rm omits \ } \Gamma\}$;

\item  if $M'$ is a $\tau_1$-substructure of $N'$ where $M', N'$ satisfy
$T_1$ and omit $\Gamma$ then $M'\restrict \tau \subm N'\restrict
\tau$; and

\item if $M \prec N \in \bK$ and $M' \in EC(T_1,
    \Gamma)$ such that $M' \restrict \tau = M$, then there is $N'
    \in EC(T_1, \Gamma)$ such that $M' \subset N'$ and $N'
    \restrict \tau = N$.

\end{enumerate}
    \end{fact}

The exact assertion  for part 3 is new in this paper; we don't
include the slight modification in the standard proofs (e.g.
\cite[Theorem 4.15]{Baldwincatmon}).  Note that we have a weakening
of Definition~\ref{funcexp} caused by the possibility of multiple
`good' expansion of a model $M$.

Here are the syntactic conditions equivalent to DAP and DJEP.

\begin{definition}
\begin{itemize}
	\item $\Psi$ has \emph{$<\lambda$-DAP satisfiability} iff for
any expansion by constants $\bc$ and all sets of atomic and
negated atomic formulas (in $\tau(\Psi) \cup \{\bc\}$)
$\delta_1(\bx, \bc)$ and $\delta_2(\by, \bc)$ of size $<
\lambda$, if $\Psi \wedge \exists \bx \left(\bigwedge
\delta_1(\bx, \bc) \wedge \bigwedge x_i \neq c_j\right)$ and
$\Psi \wedge \exists \by \left(\bigwedge \delta_2(\by, \bc)
\wedge \bigwedge y_i \neq c_j\right) $ are separately
satisfiable, then so is
	$$\Psi \wedge \exists\bx, \by\left( \bigwedge \delta_1(\bx, \bc) \wedge \bigwedge \delta_2(\by, \bc) \wedge \bigwedge_{i, j} x_i \neq y_j\right)$$

	\item $\Psi$ has \emph{$<\lambda$-DJEP satisfiability} iff
for all sets of atomic and negated atomic formulas (in
$\tau(\Psi)$) $\delta_1(\bx)$ and $\delta_2(\by)$ of size $<
\lambda$, if $\Psi \wedge \exists \bx \bigwedge \delta_1(\bx)$
and $\Psi \wedge \exists \by \bigwedge \delta_2(\by)$ are
separately satisfiable, then so is
	$$\Psi \wedge \exists\bx, \by\left( \bigwedge \delta_1(\bx) \wedge \bigwedge \delta_2(\by) \wedge \bigwedge_{i, j} x_i \neq y_j\right)$$
\end{itemize}
\end{definition}

We now outline the argument for $DJEP$; the others are similar.  Note
that $(2) \to (1)$ for the analogous result with DAP replacing DJEP has been shown by Hyttinen and
Kes{\"a}l{\"a}           \cite[2.16]{HK}.

\begin{lemma} \label{syntrans}
Suppose that $\bK$ is an AEC, $\lambda > LS(\bK)$, and $T_1$ and $\Gamma$ are from Shelah's Presentation Theorem.
  Let $\Phi$ be the $L_{LS(\bK)^+, \omega}$ theory that asserts the satisfaction of $T_1$ and omission of each type in $\Gamma$.  Then the following are equivalent:
\begin{enumerate}
	\item $\bK_{<\lambda}$ has $DJEP$.
	\item $\left(EC(T_1, \Gamma), \subset\right)_{<\lambda}$ has DJEP.
	\item $\Phi$ has $<\lambda$-$DJEP$-satisfiability.
\end{enumerate}
\end{lemma}

{\bf Proof:} \begin{enumerate}
	\item[$(1) \leftrightarrow (2)$:] First suppose that
$\bK_{<\lambda}$ has DJEP.  Let $M_0^*, M_1^* \in EC(T_1,
\Gamma)_{<\lambda}$ and set $M_\ell := M_\ell^* \rest \tau$.  By
disjoint embedding for $\ell  = 0,1$, there is $N \in \bK$ such that
each $M_\ell \prec N$. Our goal is to expand $N$ to be a member of
$EC(T_1, \Gamma)$ in a way that respects the already existing
expansions.

Recall from the proof of Fact \ref{pres} that expansions of $M \in \bK$ to models $M^* \in EC(T_1, \Gamma)$ exactly come from writing $M$ as a directed union of $LS(\bK)$-sized models indexed by $P_\omega |M|$, and then enumerating the models in the union.  Thus, the expansion of $M_\ell$ to $M^*_\ell$ come from $\{ M_{\ell, \ba} \in \bK_{LS(\bK)} \mid \ba \in M_\ell\}$, where $|M_{\ell, \ba}| = \{ \left(F_{|\ba|}^i\right)^{M^*_\ell}(\ba) \mid i < LS(\bK)\}$ and the functions $F_{n}^i$ are from the expansion.  Because $M_1$ and $M_2$ are disjoint strong submodels of $N$, we can write $N$ as a directed union of $\{N_{\ba} \in \bK_{LS(\bK)} \mid \ba \in N\}$ such that $\ba \in M_\ell$ implies that $M_{\ell, \ba} = N_{\ba}$.  Now, any enumeration of the universes of these models of order type $LS(\bK)$ will give rise to an expansion of $N$ to $N^* \in EC(T_1, \Gamma)$ by setting $\left(F_{|\ba|}^i\right)^{N^*}(\ba)$ to be the $i$th element of $|N_{\ba}|$.

Thus, choose an enumeration of them that agrees with the original enumerations from $M^*_\ell$; that is, if $\ba \in M_\ell$, then the $i$th element of $|N_{\ba}|=|M_{\ell, \ba}|$ is $\left(F_{|\ba|}^i\right)^{M^*_\ell}(\ba)$ (note that, as used before, the disjointness guarantees that there is at most one $\ell$ satisfying this).  In other words, our expansion $N^*$ will have
$$\ba \in M_\ell \implies \left(F_{|\ba|}^i\right)^{M^*_\ell}(\ba) = \left(F_{|\ba|}^i\right)^{N^*}(\ba) \text{ for all }i < LS(\bK)$$
This precisely means that $M^*_\ell \subset N^*$, as desired.  Furthermore, we have constructed the expansion so $N^* \in EC(T_1, \Gamma)$.  Thus, $\left(EC(T_1, \Gamma\right), \subset)_{<\lambda}$ has DJEP.\\

Second, suppose that $EC(T_1, \Gamma)$ has $\lambda$-DJEP.  Let $M_0
, M_1 \in \bK$; WLOG, $M_0 \cap M_1 = \emptyset$.
 Using Shelah's Presentation Theorem,
 we can expand to $M_0^*, M_1^* \in EC(T_1, \Gamma)$.  Then we can use disjoint embedding to find $N^* \in EC(T_1, \Gamma)$
 such that $M_1^*, M_2^* \subset N^*$.  By Shelah's Presentation Theorem \ref{pres}.(1), $N := N^* \rest \tau$ is the desired model.

	\item[$(2) \leftrightarrow (3)$:] First, suppose that $\Phi$ has
$<\lambda$-DJEP satisfiability.  Let $M^*_0, M^*_1 \in EC(T_1, \Gamma)$ be of size $<\lambda$.  Let
$\delta_0(\bx)$ be the quantifier-free diagram of $M_0^*$and $\delta_1(\by)$be the quantifier-free diagram of $M_1^*$.  Then $M^*_0 \vDash \Phi \wedge \exists\bx \bigwedge
\delta_0(\bx)$; similarly, $\Phi \wedge \exists \by \bigwedge \delta_1(\by)$ is satisfiable.  By the satisfiability property, there
is $N^*$ such that
$$ N^* \vDash\Psi \wedge \exists\bx, \by\left( \bigwedge \delta_0(\bx) \wedge \bigwedge \delta_1(\by) \wedge \bigwedge_{i, j} x_i \neq y_j\right)$$
Then $N^* \in EC(T_1, \Gamma)$ and contains disjoint copies of $M_0^*$ and $M_1^*$, represented by the witnesses of $\bx$ and $\by$, respectively.

Second, suppose that $(EC(T_1, \Gamma), \subset)_{<\lambda}$ has
DJEP.  Let $\Phi \wedge \exists \bx \bigwedge \delta_1(\bx)$ and
$\Phi \wedge \exists \by \bigwedge \delta_2(\by)$ be as in the hypothesis of
$<\lambda$-DJEP satisfiability.  Let $M_0^*$ witness the satisfiability of the first and $M_1^*$ witness the satisfiability of the second; note both of these are in $EC(T_1, \Gamma)$.  By DJEP, there is $N
\in EC(T_1, \Gamma)$ that contains both as substructures.  This witnesses 	$$\Psi \wedge \exists\bx, \by\left( \bigwedge \delta_1(\bx) \wedge \bigwedge \delta_2(\by) \wedge \bigwedge_{i, j} x_i \neq y_j\right)$$
Note that the formulas in $\delta_1$ and $\delta_2$ transfer up because they are atomic or negated atomic.
\end{enumerate} \hfill \dag\\

The following is a simple use of the syntactic characterization of strongly compact cardinals.

\begin{lemma} \label{sctrans}
Assume $\kappa$ is strongly compact and let $\Psi \in L_{\kappa,
\omega}(\tau_1)$ and $\lambda > \kappa$.  If $\Psi$ has
$<\kappa$-DJEP-satisfiability, then $\Psi$ has
$<\lambda$-DJEP-satisfiability.
\end{lemma}

{\bf Proof:} $<\lambda$-DJEP satisfiability hinges on the consistency of a particular $L_{\kappa,\omega}$ theory.  If $\Psi$ has $<\kappa$-DJEP-satisfiability, then every $<\kappa$ sized subtheory is consistent, which implies the entire theory is by the syntactic version  of strong compactness
we introduced at the beginning of this section.

\hfill \dag\\

Obviously the converse (for $\Psi \in L_{\infty,\omega}$) holds
without any large cardinals.

{\bf Proof of Theorem \ref{mt} for $DAP$ and $DJEP$:}  We first
complete the proof for DJEP. By Lemma \ref{syntrans}, $<\kappa$-DJEP
implies that $\Phi$ has $<\kappa$-DJEP satisfiability.  By Lemma
\ref{sctrans}, $\Phi$ has $<\lambda$-DJEP satisfiability for every
$\lambda \geq \kappa$.  Thus, by Lemma \ref{syntrans} again, $\bK$
has DJEP. The proof for DAP is exactly analogous.  \dag

\subsection{The relational presentation theorem}\label{relpres}

We modify Shelah's Presentation Theorem by eliminating the two
instances where an arbitrary choice must be made: the choice of
models in the cover and the choice of an enumeration of each covering
model. Thus the new expansion is functorial
(Definition \ref{funcexp}). However, there is a price to pay for this
canonicity. In order to remove the choices, we must add predicates of
arity $LS(K)$ and the relevant theory must allow $LS(K)$-ary
quantification, potentially putting it in $\mathbb{L}_{(2^\kappa)^+,
\kappa^+}$, where $\kappa = LS(K)$;
 contrast this with a theory of size $\leq 2^{\kappa}$ in $\mathbb{L}_{\kappa^+, \omega}$ for Shelah's version.
 As a possible silver lining, these arities can actually be brought down to $\mathbb{L}_{(I(\bK, \kappa) + \kappa)^+, \kappa^+}$.
Thus, properties of the AEC, such as the number of models in the
L{\"o}wenheim-Skolem cardinal are reflected in the presentation,
while this has no effect on the Shelah version.

We fix some notation.  Let $\bK$ be an AEC in a vocabulary $\tau$ and
let $\kappa =  LS(\bK)$.  We assume that $\bK$ contains no models of
size $< \kappa$. The same arguments could be done with $\kappa >
LS(\bK)$, but this case reduces to applying our result to $\bK_{\geq
\kappa}$.

We fix a collection of compatible enumerations for models $M \in
K_\kappa$. {\em Compatible enumerations} means that each $M$ has an
enumeration of its universe, denoted $\bm^M =\langle m_i^M : i <
\kappa\rangle$, and, if $M \cong M'$, there is some fixed isomorphism
$f_{M, M'} : M \cong M'$ such that $f_{M, M'}(m_i^M) = m_i^{M'}$ and
if $M \cong M' \cong M''$, then $f_{M, M''} = f_{M', M''} \circ f_{M,
M'}$.

For each isomorphism type $[M]_{\cong}$ and $[M \prec N]_{\cong}$ with $M, N \in K_\kappa$, we add to $\tau$
\begin{center}
$R_{[M]}(\bx)$ and $R_{[M \prec N]}(\bx; \by)$
\end{center}
as $\kappa$-ary and $\kappa 2$-ary predicates to form $\tau^*$.

A skeptical reader might protest that we have made many arbitrary choices so soon after singing the praises of our choiceless method.
 The difference is that all choices are made prior to defining the presentation theory, $T^*$.

  Once $T^*$ is defined, no other choices are made.

The goal of the theory $T^*$ is to recognize every strong submodel of
size $\kappa$ and every strong submodel relation between them via our
predicates. This is done by expressing in the axioms below concerning
sequences $\bx$ of length at most $\kappa$ the following properties
connecting the canonical enumerations with structures in $\bK$.

\begin{center}
$R_{[M]}(\bx)$ holds iff $x_i \mapsto m_i^M$ is an isomorphism

$R_{[M\prec N]}(\bx, \by)$ holds iff $x_i \mapsto m_i^M$ and $y_i \mapsto m_i^N$ are isomorphisms and $x_i = y_j$ iff $m_i^M = m_j^N$
\end{center}

Note that, by the coherence of the isomorphisms, the choice of representative from $[M]_{\cong}$ doesn't matter.  Also, we might have $M \cong M'$; $N \cong N'$; $M \prec N$ and $M' \prec N'$; but not $(M, N) \cong (M', N')$.  In this case $R_{[M \prec N]}$ and $R_{[M' \prec N']}$ are different predicates.

We now write the axioms for $T^*$. {\em A priori} they are in the
logic $\mathbb{L}_{(2^\kappa)^+, \kappa^+}(\tau^*)$ but the theorem
states a slightly finer result. To aid in understanding, we include a
description prior to the formal statement of each property.

\begin{definition} \label{tstar-def}\label{expth}
The theory $T^*$ in $\mathbb{L}_{(I(\bK, \kappa)+\kappa)^+, \kappa^+}(\tau^*)$ is the collection of the following schema:
\begin{enumerate}
	\item \label{1} {\em If $R_{[M]}(\bx)$ holds, then $x_i \mapsto m_i^M$ should be an isomorphism.}\\
	If $\phi(z_1, \dots, z_n)$ is an atomic or negated atomic $\tau$-formula that holds of $m_{i_1}^M, \dots, m_{i_n}^M$, then include
	$$\forall \bx \left( R_{[M]}(\bx) \to \phi(x_{i_1}, \dots, x_{i_n}) \right)$$

	\item \label{2} {\em If $R_{[M\prec N]}(\bx, \by)$ holds, then $x_i \mapsto m_i^M$ and $y_i \mapsto m_i^N$ should be isomorphisms and the correct overlap should occur.}\\
	If $M \prec N$ and $i \mapsto j_i$ is the function such that $m_i^M = m_{j_i}^N$, then include
	$$\forall \bx, \by \left( R_{[M \prec N]}(\bx, \by) \to
\left(R_{[M]}(\bx) \wedge R_{[N]}(\by) \wedge \bigwedge_{i <
\kappa} x_i = y_{j_i} \right) \right)$$

	\item \label{7} {\em Every $\kappa$-tuple is covered by a
model.}\\  Include the following where $\lg(\bx) = \lg(\by) =\kappa$
	$$\forall \bx \exists \by \left( \bigvee_{[M]_{\cong} \in
K_\kappa / \cong} R_{[M]}(\by) \wedge \bigwedge_{i < \kappa}
\bigvee_{j <\kappa} x_i = y_{j_i}\right)$$

	\item \label{8} {\em If $R_{[N]}(\bx)$ holds and $M\prec N$, then $R_{[M \prec N]}(\bx^{\circ}, \bx)$ should hold for the appropriate subtuple $\bx^{\circ}$ of $\bx$.}\\
	If $M \prec N$ and $\pi:\kappa \to \kappa$ is the unique map so $m_i^M = m_{\pi(i)}^N$, then denote $\bx^\pi$ to be the subtuple of $\bx$ such that $x_i^\pi = x_{\pi(i)}$ and include
	$$\forall\bx \left( R_{[N]}(\bx) \to R_{[M\prec N]}(\bx^\pi, \bx) \right)$$
	
	\item \label{9} {\em Coherence: If $M \subset N$ are both strong substructures of the whole model, then $M \prec N$.}\\
	If $M \prec N$ and $m^M_i = m^N_{j_i}$, then include
	$$\forall \bx, \by \left( R_{[M]}(\bx) \wedge R_{[N]}(\by) \wedge \bigwedge_{i < \kappa} x_i = y_{j_i} \to R_{[M \prec N]}(\bx, \by) \right)$$
\end{enumerate}
\end{definition}

\begin{remark}
We have intentionally omitted the converse to Definition \ref{tstar-def}.(\ref{1}), namely
$$\forall \bx \left( \bigwedge_{\phi(z_{i_1}, \dots, z_{i_n}) \in tp_{qf}(M/\emptyset)} \phi(x_{i_1}, \dots, x_{i_n}) \to R_{[M]}(\bx) \right)$$
because it is not true.  The ``toy example'' of a nonfinitary AEC--the $L(Q)$-theory of an equivalence relation where each equivalence class is countable--gives a counter-example.
\end{remark}

For any $M^{*} \vDash T^*$, denote $M^{*}\! \restriction \! \tau$
by $M$.

\begin{theorem}[Relational Presentation Theorem] \label{preb}
\begin{enumerate}
	\item \label{b} If $M^{*} \vDash T^*$ then $M^{*}\!
\restriction \! \tau \in K$. Further, for all $M_0 \in K_\kappa$, we
have $M^{*} \vDash R_{[M_0]}(\bm)$ implies that $\bm$
enumerates a strong substructure of  $M$.
	\item \label{a} Every $M \in K$ has a unique expansion
$M^{*}$ that models $T^*$.
	\item \label{c} If $M \prec N$, then $M^{*} \subset N^{*}$.
	\item \label{d} If $M^{*} \subset N^{*}$ both model $T^*$,
then $M \prec N$.
	\item \label{e} If $M \prec N$ and $M^* \vDash T$ such that
$M^* \rest \tau = M$, then there is $N^* \vDash T$ such that $M^*
\subset N^*$ and $N^* \rest \tau = N$.
\end{enumerate}

 Moreover, this is a functorial expansion
in the sense of Vasey \cite[Definition
3.1]{vasinfstab} and $(\Mod T^*, \subset)$ is an AEC except that it
allows $\kappa$-ary relations.
\end{theorem}

Note that although the vocabulary $\tau^*$ is $\kappa$-ary, the
structure of objects and embeddings from $(\Mod T^*, \subset)$ still
satisfies all of the category theoretic conditions on AECs, as
developed by Lieberman and Rosicky \cite{lieros}.  This is because
$(\Mod T^*, \subset)$ is equivalent to an AEC, namely $\bK$, via the
forgetful functor.

{\bf Proof:} (\ref{b}): We will build a $\prec$-directed system $\{M_{\ba} \subset M : \ba \in {}^{<\omega} M\}$ that are
members of $K_\kappa$.  We don't (and can't) require in advance that $M_{\ba} \prec M$, but this will follow from our argument.\\

For singletons $a \in M$, taking $\bx$ to be $\seq{a:i < \kappa}$ in (\ref{expth}.\ref{7}), implies that there is $M'_a
\in K_\kappa$ and $ \bm^a \in {}^{\kappa} M$ with $a \in \bm^a$ such that $M \vDash R_{[M'_a]}(\bm^a)$.  By (\ref{1}), this means that $m^a_i \mapsto m_i^{M_a'}$ is an isomorphism.  Set $M_a := \bm^a$.\footnote{We mean that we set $M_a$ to be $\tau$-structure with universe the range of $\bm^a$ and functions and relations inherited from $M_a'$ via the map above.}\\

Suppose $\ba$ is a finite sequence in $ M$ and $M_{\ba'}$ is defined
for every $\ba' \subsetneq \ba$.  Using the union of the universes as
the $\bx$ in (\ref{expth}.\ref{7}), there is some $ N
\in K_\kappa$ and $\bm^{\ba} \in {}^\kappa M$ such that
\begin{itemize}
	\item $|M_{\ba'}| \subset \bm^{\ba}$ for each $\ba' \subsetneq \ba$.
	\item $M \vDash R_{[N]}(\bm^{\ba})$.
\end{itemize}
By (\ref{expth}.\ref{8}), this means that $M\vDash R_{\bar M_{\ba'} \prec \bar
N}(\bm^{\ba'}, \bm^{\ba})$, after some permutation of the parameters.
By (\ref{2}) and (\ref{1}), this means that $M_{\ba'} \prec N$; set $M_{\ba} := \bm^{\ba}$.\\

Now that we have finished the construction, we are done.  AECs are closed under directed unions, so
$\cup_{\ba \in M} M_{\ba} \in K$.  But this model has the same universe as $M$ and is a substructure of $M$; thus
$M = \cup_{\ba \in M} M_{\ba} \in K$.  \\

For the further claim, suppose $M^* \vDash R_{[M_0]}(\bm)$.  We can redo the same proof as above with the following change: whenever $\ba \in M$ is a finite sequence such that $\ba \subset \bm$, then set $\bm^{\ba} = \bm$ directly, rather than appealing to (\ref{expth}.\ref{7}) abstractly.  Note that $\bm$ witnesses the existential in that axiom, so the rest of the construction can proceed without change.  At the end, we have
$$\bm = M_{\ba} \prec \bigcup_{\ba' \in {}^{<\omega} M} M_{\ba'} = M$$

(\ref{a}): First, it's clear that $M \in K$ has an expansion; for each $M_0 \prec M$ of size $\kappa$, make $R_{[M_0]}(\seq{m_i^{M_0}:i < \kappa})$ hold and, for each $M_0 \prec N_0 \prec M$ of size $\kappa$, make $R_{[M_0 \prec N_0]}(\seq{m_i^{M_0}:i < \kappa}, \seq{m_i^{N_0}:i < \kappa})$ hold.
  Now we want to show this expansion is the unique one.\\
Suppose $M^+ \vDash T^*$ is an expansion of $M$.  We want to show
this is in fact the expansion described in the above paragraph.  Let $M_0 \prec M$.  By (\ref{expth}.\ref{7}) and (\ref{b}) of this theorem, there is $N_0 \prec M$ and $\bn \in {}^\kappa M$ such that
\begin{itemize}
	\item $M^+ \vDash R_{[N_0]}(\bn)$
	\item $|M_0| \subset \bn$
\end{itemize}

 By coherence, $M_0 \prec \bn$.
Since $n_i \mapsto m_i^{N_0}$ is an isomorphism, there is $M^*_0
\cong M_0$ such that $M_0^* \prec N_0$.
Note that $ T^* \models \forall \bx R_{[M_0^*]}(\bx)\leftrightarrow R_{[M_0]}(\bx)$. By (\ref{expth}.\ref{8}),
$$M^+ \vDash R_{[M_0^* \prec N_0]}(\seq{m_i^{M_0}:i<\kappa}, \bn)$$
By (\ref{expth}.\ref{2}), $M^+ \vDash R_{[M_0^*]}(\seq{m_i^{M_0}:i<\kappa})$,
 which gives us the conclusion by the further part of (\ref{b}) of this theorem.

Similarly, if $M_0 \prec N_0 \prec M$, it follows that
$$M^+ \vDash R_{[M_0 \prec N_0]}(\seq{m_i^{M_0} : i < \kappa}, \seq{m_i^{N_0}:i < \kappa})$$
Thus, this arbitrary expansion is actually the intended one.

(\ref{c}): Apply the uniqueness of the expansion and the transitivity of $\prec$.\\

(\ref{d}): As in the proof of (\ref{1}), we can build
$\prec$-directed systems $\{M_{\ba} \mcolon \ba \in {}^{<\omega} M\}$
and $\{ N_{\bb}\mcolon \bb \in {}^{<\omega} N\}$  of submodels of $M$
and $N$,   so that $M_{\ba} = N_{\ba}$ when $\ba \in {}^{<\omega} M$.
From the union axioms of AECs, we see that $M \prec N$.

(\ref{e}): This follows from (\ref{c}), (\ref{d}) of this theorem and the uniqueness of the expansion.

Recall that the map $M^* \in \Mod T^*$ to $M^* \rest \tau \in \bK$ is a an abstract Morleyization  if it is a bijection such that every isomorphism $f: M \cong N$ in $\bK$ lifts to $f:M^* \cong N^*$ and $M \prec N$ implies $M^* \subset N^*$.  We have shown that this is true of our expansion.  \hfill \dag\\

\begin{remark}{\em
The use of infinitary quantification might remind the reader of the
work on the interaction between AECs and $\mathbb{L}_{\infty,
\kappa^+}$ by Shelah \cite[Chapter IV]{Shaecbook} and Kueker
\cite{Kuekeraec} (see also Boney and Vasey \cite{bvcatinflog} for
more in this area).  The main difference is that, in working with
$\mathbb{L}_{\infty, \kappa^+}$, those authors make use of the
semantic properties of equivalence (back and forth systems and
games). In contrast, particularly in the following transfer result we
look at the syntax of $\mathbb{L}_{(2^\kappa)^+, \kappa^+}$.}
\end{remark}

The functoriality of this presentation theorem allows us to give a syntactic proof of the amalgamation, etc. transfer results without assuming disjointness (although the results about disjointness follow similarly).  We focus on amalgamation and give the details only in this case, but indicate how things are changed for other properties.

Proposition~\ref{abstracttrans} applied to this context yields the following result.
\begin{prop}
$(K, \prec)$ has $\lambda$-amalgamation [joint embedding, etc.] iff
$(\Mod T^*, \subset)$ has $\lambda$-amalgamation [joint
embedding, etc.].
\end{prop}

Now we show the transfer of amalgamation between different
cardinalities using the technology of this section.

\begin{notation} \label{ad-notation}
Fix an AEC $\bK$ and the language $\tau^*$ from Theorem \ref{preb}.
\begin{enumerate}
	\item Given $\tau^*$-structures $M^*_0 \subset M^*_1, M^*_2$, we define
the {\em amalgamation diagram} $AD(M^*_1, M^*_2/M_0^*)$ to be
\begin{center}
$\{ \phi(\bc_{\bm_0}, \bc_{\bm_1})) : \phi \text{ is quantifier-free from $\tau^*$ and for  }\ell = 0 \text{ or }
1,$\\$ M^*_\ell \vDash \phi(\bc_{\bm_0},\bc_{\bm_1}), \text{ with } \bm_0 \in M^*_0 \text{ and } \bm_1 \in M^*_\ell
\}$
\end{center} in the vocabulary $\tau^* \cup \{ c_m : m \in M^*_1 \cup
M^*_2\}$ where each constant is distinct except for the common
submodel $M_0$ and
$\bc_{\bm}$ denotes the finite sequence of constants $c_{m_1}, \dots, c_{m_n}$.\\
The {\em disjoint amalgamation diagram} $DAD(M^*_1, M^*_2/M_0^*)$ is
$$AD(M^*_1, M^*_2/M_0^*) \cup \{c_{m_1} \neq c_{m_2} : m_\ell \in M^*_\ell - M^*_0\}$$

	\item Given $\tau^*$-structures $M^*_0, M^*_1$, we define the {\em joint embedding diagram} $JD(M^*_0, M^*_1)$ to be
	$$\{ \phi(\bc_{\bm})
) : \phi \text{ is quantifier-free from $\tau^*$ and for  }\ell = 0 \text{ or }
1, M^*_\ell \vDash \phi(\bc_{\bm}) \text{ with } \bm \in M^*_\ell\}$$ in
the vocabulary $\tau^* \cup \{ c_m : m \in M^*_1 \cup M^*_2\}$ where each constant is distinct.\\
The {\em disjoint amalgamation diagram} $DJD(M^*_0, M^*_1)$ is
$$AD(M^*_1, M^*_2/M_0^*) \cup \{c_{m_1} \neq c_{m_2} : m_\ell \in M^*_\ell - M^*_0\}$$
\end{enumerate}
\end{notation}

The use of this notation is obvious.

\begin{claim}\label{synap}
Any amalgam of $M_1$ and $M_2$ over $M_0$ is a reduct of a model of
$$T^* \cup AD(M_1^*, M_2^*/M_0^*)$$
\end{claim}

{\bf Proof:}  An amalgam of $M_0 \prec M_1, M_2$ is canonically expandable to an
amalgam of $M_0^* \subset M_1^*, M_2^*$, which is precisely
a model of $T^* \cup AD(M_1^*, M_2^*/M_0^*)$.  Conversely, a model of that theory will reduct to a member of $\bK$ with embeddings of $M_1$ and $M_2$ that fix $M_0$.\hfill \dag\\

There are similar claims for other properties.  Thus, we have connected amalgamation in $\bK$ to amalgamation in $(\Mod T^*, \subset)$ to a syntactic condition, similar to Lemma \ref{syntrans}.  Now we can use the compactness of logics in various large cardinals to transfer amalgamation between cardinals.  To do this, recall the notion of an amalgamation base.

\begin{definition} For a class of cardinals $\mathcal{F}$, we say $M \in K_{\mathcal{F}}$ is a {\em $\mathcal{F}$-amalgamation base}
($\mathcal{F}$-a.b.)     if any pair of models from $K_{\mathcal{F}}$ extending $M$ can be amalgamated over $M$.  We use the same rewriting conventions as in Definition \ref{mdef}.(1), e. g., writing $\leq \lambda$-a.b. for $[LS(K), \lambda]$-amalgamation base.
\end{definition}

We need to specify two more large cardinal properties.

\begin{definition}\begin{enumerate}
\item A cardinal $\kappa$ is weakly compact
    if it is strongly inaccessible and every set of $\kappa$
    sentence in $L_{\kappa,\kappa}$ that is $<\kappa$-satisfiable is
    satisfiable is satisfiable\footnote{At one time strong inaccessiblity was not required, but this
    is the current definition}.

\item A cardinal $\kappa$ is measurable if  there exists a
    $\kappa$-additive, non-trivial, $\{0,1\}$-valued measure on the power set of
    $\kappa$.

\item $\kappa$ is $(\delta, \lambda)$-strongly compact for
    $\delta \leq \kappa \leq \lambda$ if there is a
    $\delta$-complete, fine ultrafilter on
    $\Pscr_\kappa(\lambda)$.

        $\kappa$ is  $\lambda$-strongly compact if it is $(\kappa,\lambda)$-strongly compact.

\end{enumerate}
\end{definition}

This gives us the following results syntactically.

\begin{prop}\label{usescomp} \label{lc-prop}
Suppose $LS(K) < \kappa$.
\begin{itemize}
	\item Let $\kappa$ be weakly compact and $M \in K_\kappa$.  If $M$ can be written as an increasing union  $\cup_{i<\kappa} M_i$ with each $M_i \in K_{<\kappa}$ being a $<\kappa$-a.b., then $M$ is a $\kappa$-a.b.
	
	\item Let $\kappa$ be measurable and $M \in K$.  If $M$ can be written as an increasing union $\cup_{i < \kappa} M_i$ with each $M_i$ being a $\lambda_i$-a.b., then $M$ is a $\left(\sup_{i<\kappa}\lambda_i\right)$-a.b.
	
	\item Let $\kappa$ be $\lambda$-strongly compact and $M \in K$.  If $M$ can be written as a directed union $\cup_{x \in P_{\kappa}\lambda} M_x$ with each $M_x$ being a $<\kappa$-a.b., then $M$ is a $\leq \lambda$-a.b.
\end{itemize}
\end{prop}

{\bf Proof:} The proof of the different parts are essentially the
same: take a valid amalgamation problem over $M$ and
 formulate it syntactically via Claim \ref{synap} in $\mathbb{L}_{\kappa, \kappa}(\tau^*)$.  Then use the appropriate syntactic
 compactness for the large cardinal to conclude
 the satisfiability of the appropriate theory.\\

First, suppose $\kappa$ is weakly compact and $M = \cup_{i <
\kappa} M_i \in K_\kappa$ where $M_i \in K_{<\kappa}$ is a $<\kappa$-a.b. Let
$M \prec M^1, M^2$ is an amalgamation problem from $K_\kappa$.  Find
resolutions $\seq{M^\ell_i \in K_{<\kappa} : i < \kappa}$ with $M_i
\prec M^\ell_i$ for $\ell = 1, 2$.  Then

$$T^* \cup AD(M^{1*}, M^{2*}/M^*) = \bigcup_{i < \kappa} \left( T^* \cup AD(M_i^{1*}, M_i^{2*}/M_i^*)\right)$$
and is of size $\kappa$.  Each member of the union is satisfiable (by
Claim \ref{synap} because $M_i$ is a $<\kappa$-a.b.) and of size $<
\kappa$, so $T^* \cup AD(M^{1*}, M^{2*}/M^*)$ is satisfiable.  Since
$M^1, M^2 \in K_\kappa$ were arbitrary, $M$ is a $\kappa$-a.b.\\

Second, suppose that $\kappa$ is measurable and $M = \cup_{i <
\kappa} M_i$ where $M_i$ is a $\lambda_i$-a.b. Set $\lambda = \sup_{i<\kappa} \lambda_i$ and let
$M \prec M^1, M^2$ is an amalgamation problem from $K_\lambda$.  Find resolutions $ \seq{M_i^\ell \in K : i < \kappa}$ with $M_i \prec M_i^\ell$ for $\ell = 1, 2$ and $\|M_i^\ell\| = \lambda_i$.  Then

$$T^* \cup AD(M^{1*}, M^{2*}/M^*) = \bigcup_{i < \kappa} \left( T^* \cup AD(M_i^{1*}, M_i^{2*}/M_i^*)\right)$$
Each member of the union is satisfiable because $M_i$ is a $\lambda_i$-a.b.  By the syntactic characterization of measurable cardinals (see \cite[Exercise 4.2.6]{changkeisler}), the union is satisfiable.  Thus, $M$ is $\lambda$-a.b.\\

Third, suppose that $\kappa$ is $\lambda$-strongly compact and $M =  \cup_{x \in P_{\kappa}\lambda} M_x$ with each $M_x$ being a $<\kappa$-a.b.  Let $M \prec M^1, M^2$ be an amalgamation problem from $K_\lambda$.  Find directed systems $\seq{M^\ell_x \in K_{<\kappa} \mid x P_\kappa \lambda}$ with $M_x \prec M^\ell_x$ for $\ell = 1,2$.  Then

$$T^* \cup AD(M^{1*}, M^{2*}/M^*) = \bigcup_{x \in P_\kappa \lambda} \left( T^* \cup AD(M_x^{1*}, M_x^{2*}/M_x^*)\right)$$
Every subset of the left side of size $<\kappa$ is contained in a member of the right side because $P_\kappa \lambda$ is $<\kappa$-directed, and each member of the union is consistent because each $M_x$ is an amalgamation base.  Because $\kappa$ is $\lambda$-strongly compact, this means that the entire theory is consistent.  Thus, $M$ is a $\lambda$-a.b.\hfill\dag\\

From this, we get the following corollaries computing upper bounds on the Hanf number for the $\leq \lambda$-AP.

\begin{cor}
Suppose $LS(K) < \kappa$.
\begin{itemize}
	\item If $\kappa$ is weakly compact and $K$ has $<\kappa$-AP,
then $K$ has $\leq \kappa$-AP.
	\item If $\kappa$ is measurable, $\cf \lambda = \kappa$, and
$K$ has $<\lambda$-AP, then $K$ has $\leq \lambda AP$.
	\item If $\kappa$ is $\lambda$-strongly compact and $K$ has
$< \kappa$-AP, then $K$ has $\leq \lambda$-AP.
\end{itemize}
\end{cor}

Moreover, when $\kappa$ is strongly compact, we can imitate the proof
of \cite[Corollary 1.6]{MakkaiShelahcat} to show that being an
amalgamation base follows from being a
 $<\kappa$-existentially closed model of $T^*$.  This notion turns out to be the same as the notion of $<\kappa$-universally closed from \cite{Boneyct}, and so this is an alternate proof of \cite[Lemma 7.2]{Boneyct}.

\section{The Big Gap} \label{biggap-sec}

This section concerns  examples of `exotic' behavior in small
cardinalities as opposed to behavior that happens unboundedly often
or even eventually.  We discuss known work on the spectra of
existence, amalgamation of various sorts, tameness, and categoricity.

Intuitively, Hanf's principle is that if a certain property can hold
for only set-many objects then it is eventually false. He refines
this twice.  First, if  $\Kscr$ a {\em set} of  collections of
structures $\bK$ and
 $\phi_P(X,y)$ is a formula of set theory such
$\phi(\bK,\lambda)$ means some member of $\bK$ with cardinality
$\lambda$ satisfies $P$ then there is a cardinal $\kappa_P$ such that
for any $\bK\in \Kscr$, if $\phi(\bK,\kappa')$ holds for some
$\kappa' \geq \kappa_P$, then $\phi(\bK,\lambda)$ holds for arbitrarily
large $\lambda$. Secondly, he observed that if the property $P$ is
closed down for sufficiently large members of each $\bK$, then
`arbitrarily large' can be replaced by `on  a tail' (i.e.
eventually).

{\bf Existence:} Morley (plus the Shelah presentation theorem) gives
a decisive concrete example of this principle to AEC's. Any AEC in a
countable vocabulary with countable L{\"o}wenheim-Skolem number with
models up to $\beth_{\omega_1}$ has arbitrarily large models.  And
Morley \cite{Morley65a} gave easy examples showing this bound was
tight for arbitrary sentences of  $L_{\omega_1,\omega}$. But it was
almost 40 years later that Hjorth \cite{Hjorthchar, Hjorth} showed
this bound is also tight for {\em complete}-sentences of
$L_{\omega_1,\omega}$. And a fine point in his result is interesting.

We say a $\phi$ {\em characterizes} $\kappa$, if there is a model of
$\phi$ with cardinality $\kappa$ but no larger. Further, $\phi$ {\em
homogeneously} \cite{Baumgartnerchar} characterizes $\kappa$ if
$\phi$ is a complete sentence of $L_{\omega_1,\omega}$ that
characterizes $\kappa$, contains a unary predicate $U$ such that if
$M$ is the countable model of $\phi$, every permutation of $U(M)$
extends to an automorphism of $M$ (i.e. $U(M)$ is a set of absolute
indiscernibles.) and there is a model $N$ of $\phi$  with $|U(N)|
=\kappa$.

In \cite{Hjorthchar}, Hjorth found, by an inductive procedure, for
each $\alpha<\omega_1$, a countable (finite for finite $\alpha$) set
$S_\alpha$ of complete $L_{\omega_1,\omega}$-sentences such that some
$\phi_\alpha\in S_\alpha$ characterizes
$\aleph_\alpha$\footnote{Malitz \cite{malitz} (under GCH) and
Baumgartner \cite{Baumgartnerchar} had earlier characterized the
$\beth_\alpha$ for countable $\alpha$.}. This procedure was
nondeterministic in the sense that he showed one of (countably many
if $\alpha$ is infinite) sentences worked at each $\aleph_\alpha$; it
is conjectured \cite{Souldatoscharpow} that  it may be impossible to
decide in ZFC which sentence works. In \cite{BKL}, we show a
modification of the Laskowski-Shelah example (see
\cite{LaskowskiShelahatom,
 BFKL}) gives a family of $L_{\omega_1,\omega}$-sentences $\phi_r$,
such that  $\phi_r$ homogeneously  characterizes $\aleph_r$ for $r<
\omega$. Thus for the first time \cite{BKL} establishes in ZFC, the
existence of specific sentences $\phi_r$ characterizing $\aleph_r$.


{\bf Amalgamation:} In this paper, we have established a similar
upper bound for a number of amalgamation-like properties.  Moreover,
although it is not known beforehand that the classes are eventually
downward closed, that fact falls out of the proof. In  all these
cases, the known lower bounds (i. e., examples where AP holds
initially and eventually fails) are far smaller. We state the results
for countable L{\"o}wenheim-Skolem numbers, although the \cite{BKS,
KLH} results generalize to larger cardinalities.

The best lower bounds for the disjoint amalgamation property is
$\beth_{\omega_1}$ as shown in \cite{KLH} and \cite{BKS}. In
\cite{BKS}, Baldwin, Kolesnikov, and Shelah gave examples of
$L_{\omega_1,\omega}$-definable classes that had disjoint embedding
up to $\aleph_\alpha$ for every countable $\alpha$ (but did not have
arbitrarily large models). Kolesnikov and Lambie-Hanson \cite{KLH}
show that for the collection of all coloring classes (again
$L_{\omega_1,\omega}$-definable when $\alpha$ is countable)  in a
vocabulary of a fixed size $\kappa$, the Hanf number for amalgamation
(equivalently in this example disjoint amalgamation) is precisely
$\beth_{\kappa^+}$ (and many of the classes have arbitrarily large
models).  In \cite{BKL}, Baldwin, Koerwein, and Laskowski construct,
for each $r<\omega$, a {\em complete} $L_{\omega_1,\omega}$-sentence
$\phi^r$ that has disjoint $2$-amalgamation up to and including
$\aleph_{r-2} $; disjoint amalgamation and even amalgamation fail in
$\aleph_{r-1} $ but hold (trivially) in $\aleph_{r}$; there is no
model in $\aleph_{r+1}$.

The joint embedding property and the existence of maximal models are
closely connected\footnote{Note that, under joint embedding, the
existence of a maximal model is equivalent to the non-existence of
arbitrarily large models}. The main theorem of \cite{BKSoul} asserts:
 If $\langle \lambda_i: i\le \alpha <
\aleph_1\rangle$ is a strictly increasing
 sequence of characterizable cardinals  whose models satisfy JEP$(<\lambda_0)$,
there is an $L_{\omega_1,\omega}$-sentence $\psi$ such that
\begin{enumerate}
\item The models of $\psi$ satisfy JEP$(<\lambda_0)$, while JEP
    fails for all larger cardinals and AP fails in all infinite
cardinals.
\item There exist $2^{\lambda_i^+}$ non-isomorphic maximal models
    of $\psi$ in $\lambda_i^+$, for all $i\le\alpha$, but no
    maximal models in any other cardinality; and
\item $\psi$ has arbitrarily large models.
\end{enumerate}

Thus, a lower bound on the Hanf number for either maximal models of
the joint embedding property is again $\beth_{\omega_1}$. Again, the
result is considerably more complicated for complete sentences. But
\cite{BS} show that there is a sentence $\phi$  in a vocabulary with
a predicate $X$ such that if $M \models \phi$, $|M| \leq |X(M)|^+$
and for every $\kappa$ there is a model with $|M|=\kappa^+$ and
$|X(M)| = \kappa$.  Further they note that if there is a sentence
$\phi$ that homogenously characterizes $\kappa$, then there is a
sentence $\phi'$ with a new predicate $B$ such that $\phi'$ also
characterizes $\kappa$, $B$ defines a set of absolute indiscernibles
in the countable model, and there are models $M_\lambda$ for
$\lambda\leq \kappa$ such that $(|M|,|B(M_\lambda)|) =
(\kappa,\lambda)$. Combining these two with earlier results of
Souldatos \cite{Souldatoscharpow} one obtains several different ways
to show the lower bound on the Hanf number for a complete
$L_{\omega_1,\omega}$-sentence having maximal models is
$\beth_{\omega_1}$.  In contrast to \cite{BKSoul}, all of these
examples have no models beyond $\beth_{\omega_1}$.

{\bf No maximal models:}
Baldwin and Shelah \cite{BSmax} have announced
that the exact Hanf number for the non-existence of  maximal models
is the first measurable cardinal. Souldatos observed that this
implies the lower bound on the Hanf number for $\bK$ has joint
embedding of models at least $\mu$ is the first measurable.

{\bf Tameness:} Note that the definition of a Hanf number for
tameness is more complicated as tameness is fundamentally a property
of two variables: $\bK$ is $(<\chi,\mu)$-tame if for any $N \in
\bK_\mu$, if the Galois types $p$ and $q$ over $N$ are distinct,
there is an $M \prec N$ with $|M|<\chi$ and $p \restriction M \neq q
\restriction M$.

Thus, we define the \emph{Hanf number for $<\kappa$-tameness} to be the minimal $\lambda$ such that the following holds:
\begin{center}
if $\bK$ is an AEC with $LS(\bK) < \kappa$ that is $(<\kappa, \mu)$-tame for \emph{some} $\mu \geq \lambda$, then it is $(<\kappa, \mu)$-tame for arbitrarily large
 $\mu$.
\end{center}
The results of \cite{Boneyct} show that Hanf number for $<\kappa$-tameness
 is $\kappa$ when $\kappa$ is strongly compact\footnote{This can be weakened to almost
  strongly compact; see Brooke-Taylor and Rosick\'y \cite{btr} or Boney and Unger \cite{BoUn}.}.  However, this is done by showing a much stronger ``global tameness'' result that ignores the hypothesis: \emph{every} AEC $\bK$ with $LS(\bK) < \kappa$ is $(<\kappa, \mu)$-tame for all $\mu \geq \kappa$.  Boney and Unger \cite{BoUn}, building on earlier work of Shelah \cite{sh932}, have shown that this global tameness result is actually an equivalence (in the almost strongly compact form).
 Also, due to monotonicity results for tameness, the Boney results show that the Hanf number for $<\lambda$-tameness is at most the first almost strongly compact above $\lambda$ (if such a thing exists).
  The results \cite[Theorem 4.9]{BoUn} put a large restriction on the structure of the tameness spectrum for any ZFC Hanf number.  In particular, the following
  \begin{fact}
 Let $\sigma = \sigma^\omega < \kappa \leq \lambda$.    Every AEC $\bK$ with $LS(\bK) = \sigma$ is $\left(<\kappa, \sigma^{(\lambda^{<\kappa})}\right)$-tame iff $\kappa$ is $(\sigma^+, \lambda)$-strongly compact.
 \end{fact}
 This means that a ZFC (i. e., not a large cardinal) Hanf number for $<\kappa$-tameness would consistently have to avoid cardinals of the form $\sigma^{(\lambda^{<\kappa})}$ (under GCH, all cardinals are of this form except for singular cardinals and successors of singulars of cofinality less than $\kappa$).

 One could also consider a variation of a Hanf number for $<\kappa$ that requires $(<\kappa, \mu)$-tameness \emph{on a tail} of $\mu$, rather than for arbitrarily large $\mu$.  The argument above shows that that is exactly the first strongly compact above $\kappa$.
%

{\bf Categoricity:}  Another significant instance of Hanf's
observation is Shelah's proof in \cite{Sh394} that if $\Kscr$ is
taken as all AEC's $\bK$ with $LS_{\bK}$ bounded by a cardinal
$\kappa$, then there is such an eventual Hanf number for categoricity
in a successor. Boney \cite{Boneyct} places an upper bound on this
Hanf number as the first strongly compact above  $\kappa$. This
depended on  the results on tameness discussed in the previous
paragraphs.

Building on work of Shelah \cite{Shaecbook, Shaec2book}, Vasey \cite{vuniclass} proves
that if a universal class (see \cite{Sh300}) is categorical in a $\lambda$
at least the Hanf number for existence, then it has amalgamation in
all $\mu \geq \kappa$.  The he shows that for universal class in a
countable vocabulary, {\em that satisifies amalgamation}, the Hanf
number for categoricity is at most $\beth_{\beth_{(2^{\omega})^+}}$.
Note that the lower bound for the Hanf number for categoricity is
$\aleph_\omega$, (\cite{HartShelah, BKHS}).

\begin{question}
\begin{enumerate}
\item Can one calculate  in ZFC an  upper bound on these Hanf
    numbers for `amalgamation'? Can\footnote{Grossberg initiated
    this general line of research.} the gaps  in the upper and
    lower bounds of the Hanf numbers reported here be closed in
    ZFC? Will smaller large cardinal axioms suffice for some of
    the upper bounds? Does categoricity help? \item (Vasey) Are
    there any techniques for downward transfer of
    amalgamation\footnote{Note that there is an easy example in
    \cite{BKS} of a sentence in $L_{\omega_1,\omega}$ that is
    categorical and has amalgamation in every uncountable
    cardinal but it fails both in $\aleph_0$.}?
    \item Does every AEC have a functional expansion to a
        $PC\Gamma$ class.  Is there a natural class of AEC's with
        this property --- e.g. solvable groups?
\item Can\footnote{This question seems to have originated from discussions of Baldwin, Souldatos, Laskowski, and Koerwien.} one
    define in ZFC a sequence  of sentences $\phi_\alpha$ for
    $\alpha< \omega_1$,  such that
    $\phi_\alpha$ characterizes $\aleph_\alpha$?
    \item (Shelah) If $\aleph_{\omega_1}< 2^{\aleph_0}$
    $L_{\omega_1,\omega}$-sentence has models up to
    $\aleph_{\omega_1}$, must it have a model in $2^{\aleph_0}$?
    (He proves this statement is consistent in \cite{Sh522}).
    \item (Souldatos) Is any cardinal except $\aleph_0$
        characterized by a complete sentence of
        $L_{\omega_1,\omega}$ but not homogeneously?

        \end{enumerate}
    \end{question}

%

\end{document}